\documentclass[12pt]{article}

\input xy
\xyoption{all}
\usepackage{amsmath,hyperref,amssymb,amsthm,mathrsfs}
\usepackage{xcolor}
\usepackage[pdftex]{pict2e}

\usepackage[letterpaper,top=2cm,bottom=2cm,left=2.5cm,right=2.5cm,marginparwidth=1.75cm]{geometry}

\newtheorem{thm}{Theorem}[section]
\newtheorem{cor}[thm]{Corollary}
\newtheorem{lem}[thm]{Lemma}
\newtheorem{prop}[thm]{Proposition}
\newtheorem{conj}[thm]{Conjecture}
\newtheorem{ques}[thm]{Question}

\theoremstyle{remark}

\def\codim{\mathrm{codim}}

\def\E{\mathrm{E}}

\def\Sing{\mathcal{S}ing}

\def\B{\mathrm{B}}

\def\D{\mathbb{D}}
\def\Dual{\mathrm{D}}

\def\Hom{\mathrm{Hom}}

\def\F{\mathbb F}

\def\R{\mathbb R}
\def\Z{\mathbb Z}

\def\A{\mathscr A}
\def\U{\mathscr U}

\def\calR{\mathcal{R}}

\def\GL{\mathrm{GL}}

\def\ra{\rightarrow}

\def\Sq{\mathrm{Sq}}

\def\Dual{\mathrm{D}}

\begin{document}

	\title{La suspension homologique pour les CW-complexes finis quotients d'actions libres de $(\Z/2)^n$}

	\author{Nguyen Dang Ho Hai
		\and
		Lionel Schwartz
	}



		
		\maketitle

		\begin{abstract}
			
		\end{abstract}


		\section{Introduction et r\'esultats}
		
		Cette note est un compl\'ement à \cite{BLSZ}, \cite{Hai} et est \`a comparer \`a  \cite{HS}.

		Soit $G$ un groupe agissant sur un CW-complexe $X$. On notera de fa\c{c}on g\'en\'erique $\Sing$ le lieu singulier de cette action, c'est \`a dire le sous-ensemble des points dont le sous-groupe d'isotropie n'est pas r\'eduit \`a l'\'el\'ement neutre, le contexte doit rendre clair l'action consid\'er\'ee, si n\'ecessaire on notera $\Sing_G(X)$. L'action de $G$ sur $X_{reg} =X \setminus  \Sing $, la partie r\'eguli\`ere, est libre. 
		
		Soit $V_n$  un groupe ab\'elien \'el\'ementaire de rang $n$, on  abr\`egera en  $V$ dans la plupart des situations,  $W$, $U$ d\'esigneront  des sous-groupes de $V$,  $X^{W}$ d\'esignera le $V/W$-sous-espace de $X$ des points fixes sous $W$.

		On peut obtenir un mod\`ele pour le classifiant $BV$ de la mani\`ere suivante. On choisit une repr\'esentation $\calR$ de $V$, on note l'espace de repr\'esentation aussi $\calR$. On suppose   que $\calR^V=0$ et  que le lieu singulier $Sing$, qui est une r\'eunion de sous-espaces de $\calR$,  est non r\'eduit \`a $0$ (ce qui est vrai d\`es que $n \geq 2$)  et n'est pas $\calR$ tout entier. La repr\'esentation r\'eguli\`ere r\'eduite satisfait \`a ces conditions.

		Le groupe $V$ agit librement sur le lieu r\'egulier de l'action :
		$\calR_{reg}= \calR \setminus Sing$ 
		et  sur $(\calR^{\oplus k})_{reg} $, ces ensembles sont non-vides \`a  cause de la seconde hypoth\`ese sur $\calR$. Pour tout $k$ on a une inclusion  $(\calR^{\oplus k})_{reg} \hookrightarrow (\calR^{\oplus k+1})_{reg}$.
		La dimension du lieu singulier de $\calR^{\oplus k}$ croissant lin\'eairement avec $k$ et \'etant non nulle, la dualit\'e d'Alexander implique que la connexit\'e de $(\calR^{\oplus k})_{reg}$ tend vers l'infini avec $k$. La colimite est donc
		un mod\`ele pour $EV$. La colimite des quotients successifs
		$M_k=(\calR^{\oplus k})_{reg}/V$ est un mod\`ele pour $BV$.
		
		On note $c\mathcal{R}$ le compactifi\'e \`a l'infini de $\mathcal{R}$, l'action de $V$ s'y \'etend uniquement. Par projection st\'er\'eographique on a : $$c\mathcal{R} \cong S(\mathcal{R} \oplus \R) \cong  \Sigma_{nr}(S(\mathcal{R}))$$ o\`u $S(E)$ d\'esigne la sph\`ere unit\'e, on suppose que $E$ est munit d'une structure euclidienne invariante, plus g\'en\'eralement la notation vaut pour le fibr\'e en sph\`eres d'un fibr\'e vectoriel $E$ munit d'une m\'etrique, enfin $\Sigma_{nr}$ d\'esigne  la suspension non r\'eduite.
		
		La vari\'et\'e $M_k$ est  ouverte et est diff\'eomorphe \`a l'int\'erieur d'une vari\'et\'e \`a bord (\cite{BLSZ}) et le compactifi\'e \`a l'infini, $cM_k$, de $M_k$ est hom\'eomorphe au quotient de cette vari\'et\'e par son bord et au quotient 
		$(c\calR^{\oplus k}/ Sing)/V$ (voir \cite{BLSZ} section  7).

		Ces hom\'eomorphismes respectent la structure de $V_n$-espaces.
		Pour tout $k$ on a une inclusion $cM_k \hookrightarrow cM_{k+1}$. Voici un premier r\'esultat :
		
		\begin{thm} La colimite des espaces $cM_k$ est homotopiquement \'equivalente à $\Sigma^n (BV+)$.
		\end{thm}

		Dans cette note $H^*X$, $H_V^*X$ (resp. $\tilde H^*X$, $\tilde H_V^*X$) d\'esignent la cohomologie  et la cohomologie \'equivariante (resp. r\'eduite ...) modulo $2$. L'isomorphisme de Thom implique que la cohomologie $V$-\'equivariante r\'eduite 
		des espaces $cM_k$ est un un $H^*V$-module libre de rang $1$ (voir \cite{BLSZ} section 7).
		
		La question initiale, qui a motiv\'e cette note, soulev\'ee dans \cite{BLSZ}, \'etait la suivante : 	
		\begin{ques}\label{question} Soit un $V_n$-CW-complexe fini $X$ dont la cohomologie \'equivariante $H_V^*X$ est un $H^*V$-module libre. L'espace $(X/\Sing)/V_n$ est il, \`a homotopie pr\`es,  une $n$-suspension, ou r\'etracte 
			d'une $n$-suspension?
			
			Alternativement on pose aussi la question suivante : la filtration ci-dessus de $\Sigma^n (BV+)$, en un sens \`a pr\'eciser,  d\'esuspend t'elle $n$-fois?
		\end{ques}
		
		Une indication dans cette direction est que la cohomologie r\'eduite $\tilde H^*(cM_k)$ est la $n$-suspension d'un module instable  (\cite{BLSZ}). La question  est probablement beaucoup trop optimiste. Par contre on a le r\'esultat plus faible suivant :
		
		\begin{thm}\label{suspensionhomologique}  Sous l'hypoth\`ese ci dessus l'application  induite en homologie par l'\'evaluation (la suspension homologique)
			$$\Sigma^{n-1} \Omega^{n-1} (X/\Sing)/V_n \ra (X/\Sing)/V_n$$ est surjective.
		\end{thm}
		
		De ce fait on sait calculer l'homologie (et la cohomologie) de $\Omega^{n-1} (X/\Sing)/V_n$ comme foncteur de $H_*(X/\Sing)/V_n$.
		
		\vskip 5mm
		
		Dans le  cas o\`u le  $H^*V$-module libre $H_V^*X$  est libre de rang $1$, on sait, \cite{BLSZ} section 7, que $X$ est \`a homotopie pr\`es 
		le  compactifi\'e d'une repr\'esentation de $V_n$. La conjecture suivante  est alors plus raisonnable :
		
		\begin{conj} \label{suspension} 
			Soit un $V_n$-CW-complexe fini $X$ dont la cohomologie \'equivariante $H_V^*X$ est un $H^*V$-module libre de rang $1$. L'espace $(X/\Sing)/V_n$ est, \`a homotopie pr\`es,  la suspension $n$-\`eme d'un espace ou r\'etracte par d\'eformation d'une telle suspension. 
		\end{conj}
		
		Dans ce cas \ref{suspensionhomologique} peut \^etre am\'elior\'e.

		\begin{thm}\label{suspensionhomologique2} Soit $\calR$ satisfaisant aux conditions ci-dessus. On a une \'equivalence d'homotopie 
			$$
			\Sigma (	c\mathcal{R}/\Sing)/V_n  \simeq \vee_1^{d(\calR)} S^{n+1} \bigvee  \Sigma C		$$ 
			avec $C$ $n$-connexe, $d(\calR)$ un entier
			et la suspension homologique induite par l'\'evaluation 
			$$
			\Sigma^n \Omega^{n}
			C \ra  C
			$$
			est surjective.
			
		\end{thm}

		On peut encore pr\'eciser dans le cas de la repr\'esentation r\'eguli\`ere r\'eduite  $\widetilde {\calR eg}$. 
		\begin{prop}\label{suspensionhomologique3} Pour tout $d \geq 1$
			$$
			\Sigma (c\widetilde {\calR eg}^{ \oplus d}/\Sing)/V_n \simeq  (\Sigma^3 T )\vee  \Sigma X_d
			$$ 
			o\`u $T$ est $\GL(V_n)$-homotopiquement \'equivalent \`a l'immeuble de Tits,  $ X_d$ $n$-connexe et la suspension homologique induite par l'\'evaluation 
			$$\Sigma^{n} \Omega^{n}  X_d \ra  X_d$$ est surjective.
			
		\end{prop}
		
		Cela r\'esoud une question soulev\'ee dans \cite{Hai} :
		\begin{thm} Soit $e_n \in \F_2[\GL(V_n)]$ l'idempotent de Steinberg et $T(k)$ le $k$-\`eme spectre de Brown-Gitler (\cite{GLM}). 
			On a une \'equivalence d'homotopie de spectres :
			$$\Sigma^{-n} e_n.(c\widetilde{\calR eg}^{\oplus 2}/Sing)/V_n \simeq\bigvee_{i=0}^{n} T(2^i-1)$$
		\end{thm}

		La d\'emonstration de \ref{suspensionhomologique} est une application directe  de r\'esultats de \cite{BLSZ} et \cite{HK}. 
		Enfin il pourra \^etre int\'erassant de comparer \`a 
		\cite{HS} qui est tr\`es proche de cette note au moins pour ce qui concerne la section 5 et qui en est, pour part, une motivation.

		\section{La suspension homologique}

		Cette section d\'emontre le th\'eor\`eme \ref{suspensionhomologique}:
		\begin{thm} Soit $X$ un $V_n$-espace fini dont la cohomologie \'equivariante 
			$H^*_{V_n}(X)$ est un $H^*V_n$-module libre alors la suspension homologique
			$$
			H_*(	\Sigma^{n-1}\Omega^{n-1}((X/\Sing/V_n)) \ra  H_*(X/\Sing)/V_n
			$$
			est surjective.
		\end{thm}

		La d\'emonstration de ce th\'eor\`eme est une cons\'equence directe de la proposition  1.6 de Hunter-Kuhn \cite{HK}  et de \cite{BLSZ}.
		
		L'id\'ee de d\'emonstration du th\'eor\`eme \ref{suspensionhomologique} est tr\`es simple, cela en est embarassant : les complexes de \cite{BLSZ} fournissent les applications n\'ecessaires \`a l'application du r\'esultat de Hunter-Kuhn. On \'enonce cette proposition remarquable, en la rebaptisant th\'eor\`eme, les espaces et spectres sont suppos\'es de type fini pour l'homologie et $2$-complets  :
		
		\begin{thm}\label{hk}
			Soient $X$ un espace $n$-connexe, $Z$ un espace. Supposons qu'il existe une application 
			$X \ra \Sigma^nZ$ injective en homologie. Alors l'application d'\'evaluation  $\Sigma^n \Omega^n X \ra X$ est surjective en homologie.
			
			Soient $X$ un spectre $0$-connexe, $Z$ un espace. Supposons qu'il existe une application 
			$X \ra \Sigma^\infty Z$ injective en homologie. Alors l'application d'\'evaluation  $\Sigma^\infty \Omega^\infty X \ra X$ est surjective en homologie.
		\end{thm}
		
		Hunter et Kuhn  en donnent deux d\'emonstrations, l'une bas\'ee  sur la suite spectrale d'Eilenberg-Moore, l'autre sur une de Greg Arone.
		
		\medskip
		\medskip
		
		Rappelons le complexe topologique (ou complexe de Atiyah-Bredon) de \cite{BLSZ}. Celui ci \'etant exprim\'e \`a l'aide de la cohomologie \`a supports compacts on rappelle  le r\'esultat classique suivant qui permet de faire le lien (voir \cite{BLSZ} scholie 4.3) :
		\begin{prop}\label{dual}
			Soit $M$ une vari\'et\'e  \`a bord $\partial M$, $M$ compacte (et donc $\partial M$). Alors on a un isomorphisme naturel :
			$$
			H^*_c(M \setminus \partial M ) \cong H^*(M,\partial M) \cong \tilde H^*(M/ \partial M). 
			$$
		\end{prop}

		\medskip
		Soit $p$ un entier avec $-1\leq p\leq n:=\dim V$~; on d\'efinit une filtration croissante de $X$ par
		$$
		F_{p}X:=\bigcup_{\mathop{\mathrm{codim}}W\hspace{1pt}\leq\hspace{1pt}p} X^{W},
		$$
		$$
		\emptyset=F_{-1} X\subset F_{0}X=X^V  \subset \ldots\subset \ldots\subset F_{n-1} X =\Sing\subset F_nX=X.
		$$
		Le complexe de cocha\^ines ${\mathrm{C}}_{\mathrm{top}}^{\bullet}\hspace{1pt}X$ :
		$$
		\mathrm{C}_{\mathrm{top}}^{0}\hspace{1pt}X\to\mathrm{C}_{\mathrm{top}}^{1}\hspace{1pt}X\to\ldots\to\mathrm{C}_{\mathrm{top}}^{p}\hspace{1pt}X\to\ldots\to\mathrm{C}_{\mathrm{top}}^{n}\hspace{1pt}X\to 0\to \ldots
		$$
		est d\'efini par :
		$$
		\mathrm{C}_{\mathrm{top}}^{p}\hspace{1pt}X
		\hspace{4pt}:=\hspace{4pt}
		\Sigma^{-p}\hspace{2pt} H_{V}^{*}(F_{p}\hspace{1pt}X,
		F_{p-1}\hspace{1pt}X)
		\hspace{4pt}\hspace{4pt}.
		$$
		Le cobord est le connectant de la triade $(F_{p+1}X,F_{p}X,F_{p-1}X)$. 
		Ce complexe est muni d'une coaugmentation naturelle $H_{V}^*(X)=:\mathrm{C}_{\mathrm{top}}^{-1} X\to\mathrm{C}_{\mathrm{top}}^{0} X$. Le complexe coaugment\'e est not\'e $\widetilde{\mathrm{C}}_{\mathrm{top}}^{\bullet}\hspace{1pt}X$.
		Le th\'eor\`eme 0.1 de \cite{BLSZ} est le suivant~:

		\begin{thm}\label{introCtop}
			Soit  un $V$-CW-complexe fini. Les deux propri\'et\'es suivantes sont \'equivalentes~:
			\begin{itemize}
				\item[(i)] Le $H^*V$-module $H^*_VX$ est libre.
				\item[(ii)] Le complexe coaugment\'e $\widetilde{\mathrm{C}}_{\mathrm{top}}^{\bullet}\hspace{1pt}X$ est acyclique.
			\end{itemize}
		\end{thm}
		
		Il suit  alors du th\'eor\`eme 0.2 et de la proposition 0.3  de \cite{BLSZ} que :
		
		\begin{thm}
			Sous l'une des hypoyth\`eses pr\'ec\'edentes les modules sur l'alg\`ebre de Steenrod 
			$\mathrm{C}_{\mathrm{top}}^{p}X
			:=\hspace{4pt}
			\Sigma^{-p}H^*_V(F_pX,
			F_{p-1}X)$ sont instables.
		\end{thm}
		
		Pour des r\'ef\'erences sur les modules instables on renvoie \`a \cite{Sc2}.
		
		Le connectant  du complexe topologique  est induit par une application :
		$$
		EV \times_V F_{p}X/EV \times_V F_{p-1}X \ra \Sigma (EV \times_V F_{p-1}X/EV \times_V F_{p-2}X)
		$$
		o\`u les $V$-espaces $F_{k}X/ F_{k-1}X$
		ont  des point-bases $V$-stables \'evidents.
		Il  est reli\'e au connectant :
		$$
		F_{p}X/ F_{p-1}X \ra \Sigma (F_{p-1}X/F_{p-2}X)
		$$
		par le diagramme commutatif suivant, o\`u  les  colonnes des cofibrations \`a homotopie pr\`es :
		$$
		\xymatrix{
			& EV \times_V +  \ar[d]\ar[r] &  EV \times_V +   \ar[d]   	\\
			&EV \times_V F_{p}X/ F_{p-1}X \ar[r] \ar[r] \ar[d] & EV \times_V \Sigma( F_{p-1}X/ F_{p-2}X)          \ar[d]	\\
			& EV \times_VF_{p}X/ EV \times_VF_{p-1}X\ar[r]	     & \Sigma (EV \times_VF_{p-1}X/ EV \times_VF_{p-2}X ).
		}
		$$
		De plus on a un hom\'eomorphisme :
		$$
		F_{k}X/ F_{k-1}X \simeq \bigvee_{\codim W =k } X^W/Sing
		$$

		Dans le complexe topologique  la fl\`eche terminale (\`a droite) du complexe est donc induite par une application :
		$$
		EV \times _V X/\Sing  \ra  \bigvee_{\codim W =n-1 }EV \times _V \Sigma (X^W/Sing) 
		$$
		en \'ecrasant en un point de chaque c\^ot\'e le facteur $BV$ correspondant au point base. Donc \`a homotopie pr\`es on a une application:
		\begin{equation}
			(X/\Sing)/V \ra \bigvee_{\codim W = n-1} \; \Sigma(B\Z/2 \times X^W/\Sing)/(V/W)
		\end{equation}
		qui est surjective en cohomologie sous l'hypoth\`ese que le $\mathrm{H}^{*}V$-module $\mathrm{H}_{V}^{*}X$ est libre.

		Les propositions 1.1 et 4.11 de \cite{BLSZ} montrent  que si $H^*_VX$ est libre comme $H^*V$-module, alors $H^*_{V/W}X^W$ est libre comme $H^*V/W$-module.
		Donc on peut composer sur chaque facteur du bouquet ci dessus (\`a droite) par l'application obtenue en substituant, sur le facteur index\'e par $W$, $V/W$  \`a $V$  et $X^W$ \`a $X$  on obtient une application :

		\begin{equation}
			(X/\Sing)/V \ra \bigvee_{W_1 \subset W_2, \codim W_2= n-2, \,\ codim W_1 = n-1} \; \Sigma^2(B(\Z/2)^2 \times X^{W_2}/\Sing)/(V/W_2)
		\end{equation}
		qui est encore surjective en cohomologie. En it\'erant on obtient :

		\begin{lem} \label{epi}  Si la cohomologie \'equivariante de $X$ est $H^*V$-libre il existe une application 
			$$
			(X/\Sing)/V \ra \bigvee_{D} \; \Sigma^{n}(BV \times X^V +)
			$$
			vers le bouquet index\'e par les drapeaux $D=\{(0 \subsetneq W_1 \subsetneq W_2 \subsetneq \ldots \subsetneq W _n=V)$\}
			qui est surjective en cohomologie, donc injective en homologie.
			
		\end{lem}
		Ainsi qu'on l'a dit plus haut 
		\`a chaque \'etape  les conditions de libert\'e requises sont satisfaites  de par la proposition (4.11) de \cite{BLSZ}. 
		
		Notons enfin que l'espace singulier \'etant dans le dernier pas de l'it\'eration est vide, ce qui explique l'apparition du point adjoint $+$.

		On est dans les conditions d'application de la proposition 1.6 de \cite{HK}, sauf que $(X/\Sing)/V$ est seulement $(n-1)$-connexe. Ceci donne donc \ref{suspensionhomologique}. 
		
		Par contre l'application ainsi construite permettra d'appliquer directement le th\'eor\`eme 1.1 de \cite{HK}. Ceci sera pr\'ecis\'e en section 5.

		\section{Filtrations stables de $BV_n$}
		
		Dans cette section on \'etudie le cas particulier o\`u le complexe fini sur lequel agit $V_n$ est le compactifi\'e $c\calR$ d'une repr\'esentation $\calR$. On d\'emontre d'abord le th\'eor\`eme \ref{suspensionhomologique2} :
		
		\begin{thm} Soit $\calR$ une repr\'esentation de $V_n$ telle que $\calR^{V_n}=0$. On a une \'equivalence d'homotopie 
			$$
			\Sigma (	c\mathcal{R}/\Sing)/V_n  \simeq \vee_1^{d(\calR)} S^{n+1} \bigvee  \Sigma C		$$ 
			avec $C$ $n$-connexe, $d(\calR)$ un entier et la suspension homologique
			$$
			H^*(\Sigma^n \Omega^{n}
			C) \ra  H^*(C)
			$$
			est surjective.
		\end{thm}
		
		On n'a pas suppos\'e dans l'\'enonc\'e que $Sing$ est non r\'eduit \`a $0$ et non \'egal \`a $\calR$ tout entier. Dans ces deux cas l'\'enonc\'e est juste, mais trivial.
		
		Dans la suite de cette section $V_n$ sera abr\'eg\'e en $V$.
		
		On note $c\calR^k$
		le compactifi\'e de la repr\'esentation $\calR^{\oplus k}$. On a un syst\`eme inductif
		$$
		c\calR/Sing \hookrightarrow c\calR^2/Sing \hookrightarrow \ldots \hookrightarrow c\calR^k/Sing \hookrightarrow \ldots
		$$
		dont le t\'elescope sera not\'e $c\calR^\infty/Sing$. On a \'egalement le syt\`eme obtenu par quotient :
		$$
		(c\calR/Sing)/V \hookrightarrow (c\calR^2/Sing)/V \hookrightarrow \ldots \hookrightarrow (c\calR^k/Sing)/V \hookrightarrow \ldots
		$$ et le t\'el\'escope $(c\calR^\infty/Sing)/V$.
		
		\begin{lem} Si $\calR^V=0$, pour tout $k$ 
			l'application induite en cohomologie par 
			$(c\calR^k/Sing)/V \hookrightarrow (c\calR^{k+1}/Sing)/V$ est surjective, c'est un isomorphisme en degr\'e $n$.
		\end{lem}
		
		Ce lemme est d\'ej\`a d\'emontr\'e dans \cite{BLSZ} Proposition 7.39  dans le cas de la repr\'esentation r\'eguli\`ere r\'eduite. Dans le cas g\'en\'eral on la d\'emontre par r\'ecurrence sur la dimension de $V$. 
		
		On rappelle d'abord l'isomorphisme (\cite{BLSZ}) :
		$$
		{\mathrm{C}}_{\mathrm{top}}^{n}(c\calR^k) \cong \Sigma^{-n}  \tilde H^*((c\calR^k/Sing)/V)
		$$
		
		Puis on observe que si $\calR^V=0$, alors $(\calR^W)^{V/W}=0$. De plus, pour tout $n$,  l'ensemble des points fixes de l'action de $V$ sur $c\calR^n$ est r\'eduit \`a $0$ et au point \`a l'infini. 
		
		L'inclusion $\calR^{\oplus k} \hookrightarrow \calR^{\oplus k+1}$ induit un morphisme de complexes :
		$${\mathrm{C}}_{\mathrm{top}}^{\bullet}\hspace{1pt}(c\calR^{k+1}/Sing) \ra
		{\mathrm{C}}_{\mathrm{top}}^{\bullet}\hspace{1pt}(c\calR^{k}/Sing)$$
		Un argument de r\'ecurrence sur $\dim(V)$ donne alors  le r\'esultat. 
		
		En dimension $1$ on voit \`a la main que l'application est injective au niveau du terme ${\mathrm{C}}_{\mathrm{top}}^{0}$,  au niveau du terme ${\mathrm{C}}_{\mathrm{top}}^{1}$ l'application est surjecive et est l'identit\'e en degr\'e $0$. 
		
		Le r\'esultat suit de l'hypoth\`ese de r\'ecurrence en tenant compte des d\'esuspensions et \'evidemment de la comparaison des complexes topologiques pour $k$
		et $k+1$.

		\begin{prop}\label{telescope} Soit $\calR$ une repr\'esentation telle que $\calR^V=0$. Le t\'elescope $c(\calR^\infty/Sing)/V$ est homotopiquement \'equivalent 
			\`a un bouquet de $n$-suspensions de $BV+$. 
		\end{prop}

		On consid\`ere les complexes ${\mathrm{C}}_{\mathrm{top}}^{\bullet}\hspace{1pt}(c\calR^{k}/Sing)$ et on passe \`a la limite sur $k$. Le complexe obtenu  demeure exact, les conditions de Mittag-Lefler \'etant satisfaites. On peut supposer par hypoth\`ese de r\'ecurrence que tous les termes, sauf le dernier, sont des sommes directes de copies de $H^*V$. Il en est donc de m\^eme du dernier \cite{LS1}. Il suit que comme module instable, $H^*((c\calR^\infty/Sing)/V)$ est somme directe de copies de $\Sigma^nH^*V$. L'identification homotopique r\'esulte  alors des arguments de la section 2 sur le connectant du complexe.

		Ce r\'esultat donne une application 
		$$(c {\calR }^\infty/Sing)/V \ra \vee_1^{d(\calR)} S^n	$$
		qui induit un isomorphisme en cohomologie en dimension $n$ pour un certain  entier $d(\calR )$ et donc, pour $k \geq 1$ quelconque  une application :
		$$(c {\calR }^k/Sing)/V \ra \bigvee_1^{d(\calR)} S^n	$$
		qui a la m\^eme propri\'et\'e. Cette application permet de scinder le $n$-squelette
		apr\`es suspension, en notant $X_k$ le quotient de l'espace par le 
		$n$-squelette on a  :
		
		\begin{cor}\label{scindage1} Sous l'hypoth\`ese ci dessus sur $\calR$ :
			$$
			\Sigma(c{\calR }^k/Sing)/V \simeq  \vee_1^{d(\calR)} S^{n+1}  \bigvee \Sigma X_k		$$
			avec $X_k$ $n$-connexe.
		\end{cor}

		\section{L'immeuble de Tits}
		
		Cette section pr\'ecise \ref{scindage2} dans le cas de la repr\'esentation r\'eguli\`ere r\'eduite $\widetilde{\calR eg}$ et d\'emontre la proposition \ref{suspensionhomologique3}.	
		
		\begin{prop}\label{ScindageTits} Dans le cas o\`u $\calR$ est la repr\'esentation r\'eguli\`ere r\'eduite $\widetilde{\calR eg}$ on a :
			$$
			\Sigma(c\widetilde {\calR eg}^k/Sing)/V \simeq\Sigma^3 T  \vee\Sigma X_k		$$
			o\`u $T$ est $\GL (V_n)$-homotopiquement \`a l'immeuble de Tits et  $X_d$ $n$-connexe.
		\end{prop}
		
		On montre d'abord un lemme.
		
		Soit le sous-ensemble des vecteurs non nuls de $V_n$, on indexe les sommets du simplexe
		$\Delta^{2^n-1}$ par ces vecteurs. Soit $\Gamma_n$ le sous-complexe constitu\'e par les faces index\'ees par les syt\`emes de vecteurs qui n'engendrent pas $V_n$, $\Gamma_n$ est un $\GL(V_n)$-espace.
		\begin{lem}\label{Steinberg}
			Le $\GL(V_n)$-espace $\Gamma_n$ est $\GL(V_n)$-homotopiquement
			\'equivalent \`a l'immeuble de Tits de $V_n$.
		\end{lem}
		
		Ce lemme est classique, on le trouve dans des notes de Jean Lannes par exemple. En l'absence de r\'ef\'erence certaine on en donne une d\'emonstration.
		
		Soit $\mathcal S_n$ la cat\'egorie dont les objets sont les sous-ensembles de vecteurs 
		non nul  de $V_n$ qui n'engendrent pas $V_n$, les morphismes l'inclusion.
		La r\'ealisation g\'eom\'etrique du nerf, $N(\mathcal S_n)$,  de $\mathcal S_n$ est la subdivision barycentrique de $\Gamma_n$.
		
		Soit $\mathcal W_n$ la cat\'egorie des sous-espaces $W$ non-triviaux de $V_n$
		($0 \not = W \not = V_n$).
		
		Soit le foncteur $\mathcal Gen : \mathcal S_n \ra \Gamma_n$  qui envoie 
		un syt\`eme de vecteurs vers le sous-espace qu'il engendre. Il v\'erifie les conditions d'application du Th\'eor\`eme A de Quillen \cite{AKT}. En effet 
		chaque sous-cat\'egorie $W \setminus \mathcal Gen$ a pour objet terminal, l'ensemble des vecteurs non nul de $W$. Le r\'esultat suit.
		
		Rappelons le joint de $X_1, \ldots , X_t$
		$$
		\mathbb{J}(\underline{X})=\bigstar_{1 \leq i \leq t} X_i 
		$$
		qui est le quotient de :
		$$X_1 \times \ldots \times X_t \times \Delta^{t-1}$$ 
		o\`u $\Delta^{t-1}$ est le $(t-1)$-simplexe standard,  
		par la relation d'\'equivalence engendr\'ee par les relations suivantes, pour tout $i$ (o\`u ci dessous le $0$ est en position $i$) pour tous les $x_j$, $x_i,x'_i$, $a_j$,... :
		$$
		(x_1,\ldots , x_i, \ldots x_t; a_1, \ldots , 0, \ldots a_t) \sim (x_1,\ldots , x'_i, \ldots x_t; a_1, \ldots , 0, \ldots a_t) 
		$$

		Si chacun des $X_i$ est un point le joint n'est autre que $\Delta^{t-1}$. Si de plus chacun des $X_i$ est un $G$-espace, $G$ un groupe,  l'application induite par les projections sur un point :
		$$
		\mathbb{J}=\bigstar_{1 \leq i \leq t} X_i \ra \Delta^{t-1}
		$$ 
		est \'equivariante et passe  au quotient et, l'action sur $\Delta^{t-1}$ \'etant triviale, donne
		$$
		\mathbb{J}/G \ra \Delta^{t-1}
		$$ 
		Si de plus on choisit 
		pour chaque $X_i$ un point base on r\'ecup\`ere une application qui n'est pas \'equivariante :
		$$
		\Delta^{t-1} \ra \mathbb{J}
		$$
		mais par composition :
		$$
		\Delta^{t-1} \ra \mathbb{J} \ra 	\mathbb{J}/G 
		$$
		La compos\'ee
		\begin{equation}
			\Delta^{t-1} \ra \mathbb{J} \ra 	\mathbb{J}/G \ra \Delta^{t-1}
		\end{equation}
		est l'identit\'e.
		
		Explicitons dans le cas de $c{\calR eg}^k$, on abr\`egera $c{\calR eg}^k$
		en $c{\calR }$ dans ce qui suit (sauf dans le lemme). Les $2^n-1$ repr\'esentations non triviales (de dimension $1$) de $V_n$ sont not\'ees $\rho_i$, $1 \leq i \leq 2^n-1$, $\mathcal{R}$  s'\'ecrit 
		$$\mathcal{R} \cong \oplus_i \rho_i^{k}.$$
		On notera $E_i$ pour 
		$\rho_i^{\oplus k}$.
		On a  un isomorphisme de $V_n$-espaces 
		$$	
		S(\mathcal {R}) \cong \bigstar_{1 \leq i \leq 2^n-1} S(E_i) 
		$$
		o\`u $S(E_i)$ d\'esigne la sph\`ere unit\'e de $E_i$.
		
		Dans ce cas si on consid\`ere les applications consid\'er\'ees ci-dessus (en (3)) on constate que $\Gamma_{n}$ prend valeurs 
		dans $\Sing$ qui lui s'envoie sur $\Gamma_{n}$.
		L'action de $V_n$ sur $\Gamma_{n}$ \'etant triviale on obtient le :
		
		\begin{lem}\label{scindage2}
			On a un diagramme commutatif \`a homotopie pr\`es,  dont les colonnes sont de cofibrations :
			
			$$
			\xymatrix{
				& \Gamma_{n}  \ar[d]\ar[r] &   \Sing/V_n  \ar[d] \ar[r]  & \Gamma_{n}\ar[d]	\\
				&\Delta^{2^n-2}\ar[r] \ar[r] \ar[d] & c\widetilde{\calR eg}^k/V_n         \ar[d] \ar[r] & \Delta^{2^n-2} \ar[d]	\\
				& \Sigma(\Gamma_{n})\ar[r]	     & (c\widetilde{\calR eg}^k/\Sing)/V_n \ar[r] & \Sigma(\Gamma_{n})
				\\	}
			$$
			
			Les carr\'es sup\'erieurs commutent exactement. Ceux du bas  aussi
			en rempla\c{c}ant les quotients de la ligne inf\'erieure par des c\^ones.
			Les compos\'ees horizontales sup\'erieures et inf\'erieures sont l'identit\'e et $\GL(V_n)$ \'equivariante.
			
		\end{lem}
		
		La proposition \ref{ScindageTits} suit. Ceci donne une autre indication vers la conjecture \ref{suspension}. On peut en ajouter une autre. Il est clair que
		$(c\calR/\Sing)/V_n$ est une suspension, il para\^it raisonnable de montrer g\'eom\'etriquement que dans la cohomologie 
		de $\Sigma^{-1}(c\calR/\Sing)/V_n$ tous les cup-produits non triviaux sont nuls.

		\section{Les spectres de Brown-Gitler}
		
		Ce qui pr\'ec\`ede permet de d\'emontrer la conjecture \'enonc\'ee dans \cite{Hai}. Il faut faire attention \`a la terminologie ``spectres de Brown-Gitler" qui selon les sources d\'esigne le spectre ou son dual de Spanier-Whitehead. On suit ici la convention  de \cite{GLM} et de \cite{HK} -voir en particulier le commentaire au dessus de 0.4 dans \cite{GLM}- la proposition \ref{dual} permet de faire le lien. 
		Avec \cite{GLM} et \cite{HK} on appelle spectre Brown-Gitler un spectre $T(k)$ $2$-complet dont la cohomologie est isomorphe au module instable $J(k)$ et qui v\'erifie l'une des hypoth\`eses du th\'eor\`eme 1.1 de \cite{HK} dont on extrait ce qu'il nous faut ci dessous. On rappelle d'abord qu'un spectre $X$
		\`a la propri\'et\'e de Brown-Gitler si pour tous les spectres $Z$,  r\'etracte du spectre en suspension d'un espace, l'application naturelle
		$$
		[X,Z] \rightarrow{} \Hom_\A (H^*Z,H^*X)$$
		est surjective. Ici $\A$ est l'alg\`ebre de Steenrod modulo $(2)$.
		
		\begin{thm}
			Soit $X$ un spectre dont la cohomologie est injective dans la catégorie $\U$. Alors les deux conditions suivantes sont \'equivalentes.
			
			(i) Il existe un spectre $Y$ dont la cohomologie  est un module instable injectif r\'eduit et une  application   $f : X \ra Y $ surjective en cohomologie.
			
			(ii) X satisfait à la propri\'et\'e de Brown-Gitler.
		\end{thm}
		
		Le spectre que l'on va d\'ecrire arrive d'embl\'ee avec une application satisfaisant 
		aux propri\'et\'es de la premi\`ere condition ainsi qu'il est dit en fin de la section 3.
		
		Il faut montrer que sa cohomologie est celle attendue. Cela est fait dans \cite{Hai}, revenons dessus. Comme plus haut on écrit $V$ pour $V_n=(\Z/2)^n$ dans ce qui suit.	Le groupe lin\'eaire $\GL(V)$ agit sur le quotient
		$(c\widetilde{\calR eg}^2/Sing)/V$. 
		On peut, et apr\`es suspension,  appliquer l'idempotent de Steinberg et d\'ecomposer l'espace en bouquet de deux espaces.
		En fait dans ce cas il n'est pas n\'ecessaire de suspendre car 	$(c\widetilde{\calR eg}^2/Sing)/V$
		est d\'ej\`a une suspension compatible \`a l'action de $\GL(V)$.

		\begin{thm}
			Soit $e_n$ l'idempotent de Steinberg de $\F_2[\GL(V)]$. Le spectre $$\Sigma^{-n} e_n.(c\widetilde{\calR eg}^2/Sing)/V$$ est homotopiquement \'equivalent au bouquet $\bigvee_{i=0}^{n} T(2^i-1)$.
		\end{thm}
		
		Le lemme \ref{epi} dit exactement que l'on peut appliquer le th\'eor\`eme 5.1 en utilisant les r\'esultats de \cite{Hai}. 
		
		On rappelle  en quelques lignes pourquoi le spectre en question a la bonne cohomologie. 
		\'Etant donn\'e un $\A$-module $M$, on rappelle que le dual de Spanier-Whitehead $\D M$ est d\'efini par
		$$\begin{cases}
			(\D M)^{-n}=\Hom_{\F_2}(M^n,\Z/2),\qquad n\in \Z,\\
			\theta(f)=f\circ (\chi(\theta)),\qquad f\in \Dual M, \theta\in \A,
		\end{cases}$$
		$\chi$ \'etant l'antipode de $\A$. 
		On a des isomorphismes 
		\begin{eqnarray*}
			\tilde H^*((c\widetilde{\calR eg}^2/Sing)/V)&\cong& H^*_c((\widetilde{\calR eg}^2 \setminus Sing)/V) \\
			&\cong& \Sigma^{2(2^n-1)} \D H^*((\widetilde{\calR eg}^2 \setminus Sing)/V).
		\end{eqnarray*}
		Le premier isomorphisme est bien $\A$-lin\'eaire et on va expliquer la $\A$-lin\'earit\'e du second qui est donn\'e par la dualit\'e de Poincar\'e.

		Pour d\'eterminer $H^*((\widetilde{\calR eg}^2 \setminus Sing)/V)$, on note $V^*$ le dual lin\'eaire de $V$ et on identifie la repr\'esentation r\'eguli\`ere r\'eduite $\widetilde{\calR eg}$ \`a l'espace vectoriel $\R[V^*\setminus 0]$ muni de l'action de $V$ donn\'ee par 
		$$v[\alpha]=(-1)^{\alpha(v)}[\alpha], \quad v\in V, \alpha \in V^*.$$ 
		Avec cette identification, on v\'erifie que l'espace $(\widetilde{\calR eg}^2 \setminus Sing)/V$ est de la forme $\mathcal Z_K(\R,\R^*)/V$ o\`u $K$ est le complexe simplicial $\Delta(V^*,2)$ d\'efini dans \cite{Hai}, Definition 5, et $\mathcal Z_K$ le produit polyh\'edral \cite{BBC}. La cohomologie modulo $2$ de $\mathcal Z_K(\R,\R^*)/V$ dans le cas o\`u $K$ est Cohen-Macaulay a \'et\'e calcul\'ee dans \cite{DavidJanuskewiz91}. Dans notre cas, celle-ci est isomorphe \`a l'anneau quotient fini $\mathbf R(V^*,2)$ de $H^*V$ \'etudi\'e dans \cite{Hai}. On y a d\'ej\`a montr\'e que le facteur de $\mathbf R(V^*,2)$ associ\'e \`a la repr\'esentation de Steinberg est isomorphe \`a la somme directe
		$$\bigoplus_{i=0}^n\Sigma^{2(2^n-1)-n-(2^i-1)}B(2^i-1)$$
		o\`u les modules de Brown-Gitler $B(k)$ et $J(k)$ sont li\'es par $\Sigma^k \D B(k)\cong J(k)$.

		Pour finir, on explique la $\A$-lin\'earit\'e de la dualit\'e de Poincar\'e 
		$$H^*_c((\widetilde{\calR eg}^2 \setminus Sing)/V) 
		\cong\Sigma^{2(2^n-1)} \D H^*((\widetilde{\calR eg}^2 \setminus Sing)/V).$$ 
		Pour cela, soit $X$ une vari\'et\'e sans bord de dimension $d$ (pas n\'ecessairement compact, voir par exemple \cite{Hat}, Chapter 3).
		La dualit\'e de Poincar\'e est donn\'ee par $$\Dual: H_c^{d-i}X \to \Hom_{\F_2}(H^{i}X,\F_2),\qquad \Dual(x)(y)=x\cup y \in H^d_c(X)\cong \F_2.$$ 
		Pout tout $x \in H_c^{d-i}(X)$, on a $\Sq^i(x)=x\cup v_i(\tau_X)$ o\`u $v_i(\tau_X)$ est la $i$-\`eme classe de Wu du fibr\'e tangent, $\tau_X$, de $X$ (\cite{MS}). On observe que la dualit\'e de Poincar\'e $\Dual$ est $\A$-lin\'eaire si $v_i(\tau_X)=0$ pour tout $i>0$. En effet,
		pour tout $u\in H_c^{d-t-m} (X)$ et tout $v\in H^{t}(X)$, on a 
		$$\Sq^m  (u)\cup  v+\sum_{i=1}^m\Sq^i(u\cup (\chi \Sq^{m-i}v))+ u\cup \chi \Sq^m(v)=0 $$
		Donc si les classes de Wu de $\tau_X$ sont triviales le terme au milieu s'annule ce qui implique que $$\Sq^m (u)\cup v=u\cup \chi \Sq^m(v).$$
		Il suit que $\Dual$ envoie $\Sq^m(u)$ sur la fonction $$v\mapsto \Sq^m(u)\cup v=u\cup \chi\Sq^m (v),$$
		et $\Sq^m(\Dual(u))$ est la fonction compos\'ee
		$v\mapsto \chi\Sq^m (v)\mapsto u\cup \chi\Sq^m (v)$. Ceci montre que $\Dual$ commute avec $\Sq^m$.

		On revient \`a notre cas o\`u $X=(\widetilde{\calR eg}^2 \setminus Sing)/V$ est une vari\'et\'e de dimension $d=2(2^n-1)$. 
		Par instabilit\'e, si $i>\frac{d}{2}=2^n-1$ et $x\in H^{d-i}_c(X)$, alors $\Sq^i(x)=x\cup v_i(\tau_X)=0$, ce qui implique que $v_i(\tau_X)=0$ car la forme bilin\'eaire $H^{d-i}_cX\times H^iX\to H^{d}_cX\cong \F_2$ est non-d\'eg\'en\'er\'ee.

		Le fibr\'e tangent $\tau_X$ est le ``pullback" du fibr\'e $\E V\times_V \widetilde{\calR eg}^2 \xrightarrow{\rho_2} \B V$ induit par une application $f:X\to \B V$. La classe totale de Stiefel-Whitney de $\rho_2$ est donn\'ee par $$w(\rho_2)=\prod_{0\ne \alpha\in V^*} (1+\alpha)^2=(1+Q_{n,0}+\cdots +Q_{n,n-1})^2,$$ 
		o\`u $|Q_{n,i}|=2^n-2^i$ et $\F_2[Q_{n,0},\ldots,Q_{n,n-1}]\cong H^*V^{\GL_n}$ est l'alg\`ebre d'invariants de Dickson. Donc $w_i(\rho_2)$ est trivial si $0<i\le 2^{n}-1$. Par naturalit\'e, $w_i(\tau_X)=0$ si $0<i\le 2^{n}-1$. Utilisant la formule $v(\tau_X) = \chi\Sq (w(\tau_X))$, on voit que $v_i(\tau_X)=0$ si $0<i\le 2^n-1$. Donc $v_i(\tau_X)=0$ pour tout $i>0$.

		\noindent Nguyen Dang Ho Hai, \\
		College of Sciences, Hue University,\\
		77 Nguyen Hue, Hue city, VIETNAM
		\\
		Email: nguyendanghohai@husc.edu.vn
		
		\vskip 1cm
		
		\noindent Lionel Schwartz, \\
		LAGA, UMR 7539 du CNRS, Universit\'e Sorbonne Paris Nord,\\
		99 Av. J-B Cl\'ement, 93430 Villetaneuse, FRANCE \\
		IRL FVMA du CNRS. \\
		Email: schwartz@math.univ-paris13.fr
		
	\end{document}